\documentclass{article}

\usepackage[english]{babel}

\usepackage{amsmath,amssymb,setspace}
\usepackage{amsthm}
\usepackage[english]{babel}
\usepackage{amsfonts}
\usepackage{calligra}
\usepackage[T1]{fontenc}
\usepackage{xstring}

\usepackage[colorlinks,citecolor=blue,urlcolor=blue,bookmarks=false,hypertexnames=true]{hyperref}
\hypersetup{
    colorlinks=true,
    linkcolor=blue,
    filecolor=magenta,
    urlcolor=cyan,
}
\usepackage{graphics}
\usepackage[numbers,sort&compress]{natbib}
\usepackage{subcaption}
\usepackage{multicol}
\date{}
\usepackage{longtable}
\usepackage{float}
\usepackage[margin=0.7 in]{geometry}
\usepackage{graphicx}
\usepackage{mathtools}
\usepackage{amsmath}

\newtheorem{definition}{Definition}
\theoremstyle{plain}
\theoremstyle{definition}
\theoremstyle{remark}

\newtheorem{example}{Example}

\title{Fractal Nambu Mechanics: Extending Dynamics with Fractal Calculus}
\author{Alireza Khalili Golmankhaneh$^1$, Cemil Tun\c{c}$^2$, Davron Aslonqulovich Juraev $^{3,4}$\\
$^1$ Department of Physics, Urmia Branch, \\Islamic Azad University, Urmia 63896,  Iran\\
alirezakhalili2002@yahoo.co.in\\
$^2$ Department of Mathematics, Faculty of Sciences, \\ Van Yuzuncu Yil University, 65080-Campus, Van-Turkey\\
cemtunc@yahoo.com\\
$^3$ Department of Scientific Research, Innovation and Training of Scientific and Pedagogical Staff, \\ University of Economics and Pedagogy, Karshi, 180100, Uzbekistan\\
$^4$  Department of Mathematics, Anand International College of Engineering, Jaipur, 303012, India\\
juraevdavron12@gmail.com
}

\begin{document}

\maketitle

\let\thefootnote\relax
\footnotetext{ MSC2020:28A80} 
\footnote{Corresponding author Alireza Khalili Golmankhaneh}
\begin{abstract}
In this paper, we extend the principles of Nambu mechanics by incorporating fractal calculus. This extension introduces Hamiltonian and Lagrangian mechanics that incorporate fractal derivatives. By doing so, we broaden the scope of our analysis to encompass the dynamics of fractal systems, enabling us to capture their intricate and self-similar properties. This novel approach opens up new avenues for understanding and modeling complex fractal structures, thereby advancing our comprehension of these intricate phenomena. 
\end{abstract} 

\section{Introduction}
Fractal geometry explains shapes with self-similar properties that have fractional and complex numbers dimensions whose topological dimension is less than their fractal dimension \cite{Mandelbrot-1,Michel-j}. Fractals have different measures than Euclidean shapes. Fractals are seen in nature everywhere like clouds in the human body and the path of electrons and structures of materials in fractal forms show new results \cite{Barnsley-book,linton2021fractals}. Their properties have been studied with the help of mathematical analysis \cite{el1993certain}.
Fractals are explained in a mathematically rigorous manner, with an emphasis on examples and fundamental concepts such as Hausdorff dimension, self-similar sets, and Brownian motion \cite{bishop2017fractals}.
Fractal structures are seen and used in many branches of physical science and art engineering \cite{kotov2023application,bevilacqua2023inverse,heider2022self,hambly2023dimension,czachor2023contra,deppman2022tsallis,tarraf2023fractal}.
Fractal space-time was considered in the different processes to obtain more general models to explain experimental data  \cite{ord1983fractal,nottale1993fractal,Vrobel,Shlesinger-6}.
Different techniques were used to analyze fractals, including probabilistic techniques, measure theory, fractional calculus, harmonic analysis, and fractional spaces \cite{ma-5,uchaikin2013fractional,ma-6,trifce2020fractional,ma-9,Bongiorno23,jiang1998some,withers1988fundamental,bongiorno2011henstock,bongiorno2015fundamental,bongiorno2018derivatives,giona1995fractal}.
The $\alpha$-fractal function with variable parameters was investigated, focusing on the Weyl-Marchaud variable order fractional derivative under specific scaling factor conditions \cite{priyanka2023alpha,agathiyan2023integral,agathiyan2023remarks,valarmathi2023variable}. The extension of Darcy's law from conventional calculus to fractal calculus in order to quantify fluid flow in subterranean heterogeneous reservoirs was investigated \cite{damian2023mechanical,samayoa2022hausdorff}

A generalization of Hamiltonian mechanics, Multi-Hamiltonian Nambu mechanics has some applications \cite{nambu1973generalized,guha2004applications,guha2002applications}.
The flows created by a Hamiltonian over a Poisson manifold form the basis of Nambu mechanics. The trajectory in phase space is the intersection of the $N-1$ hypersurfaces specified by the invariants \cite{nambu1973generalized,nambu1952lagrangian}. The flow is perpendicular to all $N-1$ gradients of these Hamiltonians, parallel to the generalized cross product described by the corresponding Nambu bracket, and lagrangian to all $N-1$ gradients of these Hamiltonians. Nambu-Hamilton equations of motion, incorporating several Hamiltonians, describe the dynamics of Nambu mechanics. As a consistency condition, the basic identity for the Nambu bracket-which is comparable to the Jacobi identity-is presented \cite{takhtajan1994foundation,chandre2023classical,llibre2023generalized}.
The Hamiltonian embedding for Nambu's novel mechanical equations is given. Using Dirac's singular formalism, the embedding is investigated in both classical and quantization forms. The embedding Hamiltonian's transformations are contrasted with Nambu's canonical transformations \cite{bayen1975remarks,gautheron1996some}. Nambu-Poisson and Nambu-Jacobi brackets are defined inductively \cite{grabowski1999remarks}. In Nambu's generalized mechanics, the action integral is defined on a world sheet and there is an additional degree of freedom in phase space
\cite{lassig1997constrained,makhaldiani2011nambu}. The Nambu mechanics can, in theory, describe physical systems for which the Dirac mechanics is inappropriate since Nambu's mechanics is independent of Dirac's mechanics, at least in the odd instances \cite{cohen1975nambu}. The conditions for a Hamiltonian system to possess $n$-linear Nambu brackets is investigated. The system should only have one degree of freedom for $n=3$, but the criterion still holds for higher orders ($n\geq 4$). Levi-Civita tensors produce the canonical Nambu bracket, which is the only fundamental solution locally. \cite{chandre2023classical}. A geometric formulation for a Hamiltonian system-like generalized form of Nambu mechanics is provided. It introduces a nondegenerate 3-form connected to a 3n-dimensional phase space to model temporal evolution. A Poisson bracket of two forms is additionally included, enabling the description of connected rigid bodies and non-integrable fluid flow \cite{pandit1998generalized}. Nambu mechanics is discussed using numerous instances of mechanical systems with the Nambu bracket. The relationship between a mechanical system's superintegrability and the existence of the Nambu bracket is the main idea examined. The Nambu bracket's uses in field theory are presented \cite{shishanin2018nambu}. Using fractional-order differential forms and exterior derivatives, the creation of a fractional generalization of Nambu mechanics is shown. The generalized Pfaffian equations and offers a thorough analysis of one such example \cite{baleanu2009fractional,xu2015fractional,khalili2011fractional}. A fractal method to quantum chromodynamics is used to produce the Nambu-Jona-Lasino model, which includes a running coupling. Tsallis non-extensive statistical distributions in high-energy collisions are explained by the modified model \cite{megias2022nambu}.\\
The fractal analysis has been explored by multiple researchers through various approaches, including harmonic analysis \cite{ma-8}, measure theory \cite{freiberg2002harmonic,jiang1998some,giona1995fractal,bongiorno2015fundamental,bongiorno2011henstock,withers1988fundamental}, probabilistic methods \cite{ma-7}, fractional space \cite{stillinger1977axiomatic}, fractional calculus \cite{ma-6,gowrisankar2021fractal}, and non-standard methods \cite{Nottale-4}. Through an expanded framework that contains differential equations, fractal calculus expands classical calculus. These equations' solutions provide functions with fractal support, such as fractal sets and curves  \cite{parvate2009calculus,parvate2011calculus,AD-2,banchuin20224noise,Alireza-book,balankin2023vector}.
The usage of fractal versions of Newton, Lagrange, Hamilton, and Appell's mechanics is suggested. Following the definition of fractal velocity and acceleration, the Langevin equation on fractal curves is obtained. Hamilton's mechanics on fractal curves are created via the Legendre transformation, allowing for the simulation of non-conservative systems with fractional dimensions \cite{golmankhaneh2023classical}.
Nonstandard analysis incorporating hyperreal and hyperinteger numbers is used to define left and right limits and derivatives on fractal curves. The use of nonstandard analysis is used to define the fractal integral and differential forms \cite{khalili2023non}.
In this study, the generalization of Nambu mechanics is achieved by utilizing fractal derivatives and fractal differential forms on fractal manifolds. Through this approach, we extend the framework of Nambu mechanics to encompass fractal geometries and dynamics.\\
The following is the paper's outline:\\
We present a concise summary of fractal calculus in Section \ref{S-1}, outlining its key concepts and principles. Moving on to Section \ref{S-2}, we introduce fractal Nambu mechanics, explaining its fundamental ideas and framework. In Section \ref{S-3}, we demonstrate the application of fractal differential forms in deriving Nambu mechanics. Additionally, in Section \ref{S-4}, we showcase the practical implementation and usefulness of this formalism through relevant examples and applications. Finally, in Section \ref{S-5}, we provide a comprehensive conclusion that summarizes the main findings and contributions of our work in the context of fractal Nambu mechanics.
\section{Basic Tools \label{S-1}}
Here, we provide a summary of the fractal calculus for middle-$\epsilon$ Cantor set \(F \subset [c_{1}, c_{2}] \subset \mathbb{R}\) \cite{parvate2009calculus,Alireza-book}.

\begin{definition}
The indicator function of $F$ is defined as
\[
f(F, I) =
\begin{cases}
    1, & \text{if } F \cap I \neq \emptyset, \\
    0, & \text{otherwise},
\end{cases}
\]
where $I = [c_{1}, c_{2}] \subset \mathbb{R}$.
\end{definition}

\begin{definition}
The coarse-grained measure of $F \cap [c_{1}, c_{2}]$ is defined as
\begin{equation}
\mu_{\delta}^{\alpha}(F, c_{1}, c_{2}) = \inf_{|\mathcal{P}| \leq \delta} \sum_{i=0}^{n-1} \Gamma(\alpha+1)(x_{i+1} - x_i)^{\alpha} f(F, [x_i, x_{i+1}]),
\end{equation}
where $|\mathcal{P}| = \max_{0 \leq i \leq n-1} (x_{i+1} - x_i)$, and $0 < \alpha \leq 1$.
\end{definition}

\begin{definition}
The measure function of $F$ is defined as
\[
\mu^{\alpha}(F, c_{1}, c_{2}) = \lim_{\delta \rightarrow 0} \mu_{\delta}^{\alpha}(F, c_{1}, c_{2}).
\]
\end{definition}

\begin{definition}
The fractal dimension of $F \cap [c_{1}, c_{2}]$ is defined as
\begin{align*}
\dim_{\mu}(F \cap [c_{1}, c_{2}]) &= \inf\{\alpha : \mu^{\alpha}(F, c_{1}, c_{2}) = 0\} \\
&= \sup\{\alpha : \mu^{\alpha}(F, c_{1}, c_{2}) = \infty\}.
\end{align*}
\end{definition}

\begin{definition}
The integral staircase function of $F$ is defined as
\[
S_{F}^{\alpha}(x) =
\begin{cases}
    \mu^{\alpha}(F, c_{0}, x), & \text{if } x \geq c_{0}, \\
    -\mu^{\alpha}(F, x, c_{0}), & \text{otherwise},
\end{cases}
\]
where $c_{0} \in \mathbb{R}$ is a fixed number.
\end{definition}

\begin{definition}
For a function $h$ on an $\alpha$-perfect fractal set $F$, the $F^{\alpha}$-derivative of $h$ at $x$ is defined as
\[
D_{F}^{\alpha}h(x) =
\begin{cases}
    \underset{y \rightarrow x}{F_{-}\text{lim}} \frac{h(y)-h(x)}{S_{F}^{\alpha}(y)-S_{F}^{\alpha}(x)}, & \text{if } x \in F, \\
    0, & \text{otherwise},
\end{cases}
\]
if the fractal limit $F_{-}\text{lim}$ exists.
\end{definition}

\begin{definition}
The $F^{\alpha}$-integral of a bounded function $h(x)$, where $h \in B(F)$ (i.e., $h$ is a bounded function of $F$), is defined as
\begin{align*}
\int_{a}^{b} h(x) \, d_{F}^{\alpha}x &= \sup_{\mathcal{P}_{[c_{1}, c_{2}]}} \sum_{i=0}^{n-1} \inf_{x \in F \cap I} h(x) \, (S_{F}^{\alpha}(x_{i+1}) - S_{F}^{\alpha}(x_i)) \\
&= \inf_{\mathcal{P}_{[c_{1}, c_{2}]}} \sum_{i=0}^{n-1} \sup_{x \in F \cap I} h(x) \, (S_{F}^{\alpha}(x_{i+1}) - S_{F}^{\alpha}(x_i)),
\end{align*}
where $x \in F$, and the infimum or supremum is taken over all subdivisions $\mathcal{P}_{[c_{1}, c_{2}]}$.
\end{definition}

\begin{definition}
The fractal total derivative of a function $f:F^{n}\rightarrow R^{n}$ can be defined when $d_{F}^{\alpha}x_1, d_{F}^{\alpha}x_2, ..., d_{F}^{\alpha}x_n$ represent infinitesimal increments along the coordinate directions. It is denoted as $d_{F}^{\alpha} f$ and can be calculated as:
\begin{equation}
d_{F}^{\alpha} f=\sum_{i=1}^{n}(D_{F,x_{i}}^{\alpha}f) d_{F}^{\alpha}x_{i}
=\sum_{i=1}^{n}(\partial^{\alpha}_{F,x_{i}}f )d_{F}^{\alpha}x_{i}.
\end{equation}
Here, $(D_{F,x_{i}}^{\alpha}f)$ or $(\partial^{\alpha}{F,x{i}}f)$ represents the fractal partial derivative of $f$ with respect to $x_i$, and $d_{F}^{\alpha}x_{i}$ represents the infinitesimal increment in the $i$-th coordinate direction \cite{khalili2023non}.
\end{definition}
\section{Fractal Nambu mechanics \label{S-2}}
Nambu mechanics is a mathematical framework for explaining the dynamics of physical systems including several variables. Fractal Nambu mechanics is a modified form of Nambu mechanics. The idea of fractals is included in Fractal Nambu mechanics, enabling the characterization and study of complex systems with complicated and self-similar structures. The indicator function is used in Fractal Nambu mechanics to determine whether a fractal set is present or absent within a specified interval. The coarse-grained measure, which takes into account various scales of resolution, quantifies the size or extent of the fractal set within the interval. As the resolution becomes closer to zero, the measure function offers a limit that captures the fundamental characteristics of the fractal set. The complexity or self-similarity of the fractal set is gauged by the fractal dimension. The integral staircase function, which represents the total measure of the fractal set up to a specific point, is also introduced by Fractal Nambu mechanics. The fractal derivative and integral operators are defined using this function. The fractal integral computes the cumulative impact of a function on the fractal set, whereas the fractal derivative determines the rate of change of a function on the fractal set. Fractal Nambu mechanics, which incorporates fractal ideas into Nambu mechanics, offers a potent mathematical framework for investigating and comprehending complicated systems that exhibit fractal characteristics. It provides a deeper understanding of the dynamics and characteristics of such systems, allowing for a more thorough study and modeling strategy.\\
We present the mathematical framework of fractal Nambu mechanics as follows:

\begin{definition}
A fractal manifold denoted as $M^{\alpha}$, characterized by fractal dimensions, is classified as a Fractal Nambu-Poisson manifold when a mapping denoted by
\begin{equation}\label{sewaq}
\{,...,\}:[C^{\infty}(F)]^{\otimes^{n}} \mapsto C^{\infty}(F),
\end{equation}
exists. Here, $C^{\infty}(F)$ represents the space of smooth functions defined on fractal sets. This mapping can be referred to as the fractal Nambu-Poisson bracket of order $n$, defined for any $f_1, f_2, . . . , f_{2n-1}\in C^{\infty}(F)$ \cite{parvate2009calculus,Alireza-book}.

The fractal Nambu-Poisson bracket exhibits specific properties, which include the following:
\begin{itemize}
  \item Skew-symmetry:
\begin{equation}\label{wwqa}
  \{f_1, f_2, . . . , f_{n}\}=(-1)^{\epsilon(\sigma)}\{f_{\sigma(1)}, f_{\sigma(2)}, . . . , f_{\sigma(n)}\}.
\end{equation}
Here, $\sigma$ represents a permutation in the symmetric group $S_n$, which consists of all possible arrangements of $n$ elements. The function $\epsilon(\sigma)$ determines the parity of the permutation, determining whether it is even or odd. The skew-symmetry property ensures that the resulting bracket remains consistent even when the functions' order is altered, accounting for the sign change based on the parity of the permutation used.
  \item Jacobi identity:
\begin{align}\label{wqq}
 & \{\{f_{1},f_{2},f_{3},...,f_{n}\},f_{n+1},...,\nonumber\\&f_{2n-1}\}+\{f_{n},\{
f_{1},f_{2},...,f_{n+1}\},f_{n+2},...,f_{2n-1}\}=\{f_{1},...,f_{n-1},\{f_{n},...,
f_{2n-1}\}\},
\end{align}
 The fractal Nambu-Poisson bracket obeys the Jacobi identity, which is a fundamental property of Poisson brackets. It states that the bracket of the bracket of $n$ functions with $n$ additional functions, along with the bracket of the first function with the bracket of the remaining $2n-1$ functions, is equal to the bracket of the bracket of the first $n-1$ functions with the last $2n-1$ functions.
  \item  Leibniz Rule:
\begin{equation}\label{sswwqa}
  \{f_{1},f_{2},f_{3},...,f_{n+1}\}=f_{1}\{f_{2},f_{3},...,f_{n+1}\}+
f_{2}\{f_{1},f_{3},...,f_{n+1}\}.
\end{equation}
The fractal Nambu-Poisson bracket satisfies the Leibniz rule, also known as the generalized chain rule. This rule governs how the bracket operates on the product of functions, ensuring that the bracket operation can be applied to each function individually while preserving the appropriate sign and structure.
\end{itemize}
By satisfying these properties, the fractal Nambu-Poisson bracket of order $n$ provides a fundamental tool for describing the dynamics and relationships of smooth functions defined on Fractal Nambu-Poisson manifolds.
\end{definition}
\subsection{Fractal Nambu equations}
In this subsection, the fractal Nambu equations are formulated using the fractal Nambu-Poisson bracket, as shown below:
\begin{eqnarray}\label{wqqaq}
 \frac{d_{F}^{\alpha} x^{i}}{d_{F}^{\alpha}t}&=&\{x^{i},H_{1},H_{2},...,H_{n-1}\}\nonumber\\
&=& \sum_{jk...l}\epsilon_{ijk...l} D_{F,x_{j}}^{\alpha}H_{1}D_{F,x_{k}}^{\alpha}H_{2}...D_{F,x_{l}}^{\alpha}H_{n-1}\nonumber\\
&=& \frac{\partial^{\alpha}_{F}(x^{i},H_{1},H_{2},...,H_{n-1})}{\partial^{\alpha}_{F}(x^{1},x^{2},...,x^{n} )}.
\end{eqnarray}
In this formulation, the local coordinates on the fractal Nambu-Poisson manifold are represented by $x^i$, and the corresponding fractal Nambu-Hamiltonians are denoted as $H_1, H_2, . . . , H_{n-1}$. The Levi-Civita tensor, denoted as $\epsilon_{ijk...l}$, is involved in the equations.\\
It is important to note that while traditional Hamiltonian mechanics describes the state space using a pair of variables, $x$ and $p$, Nambu mechanics extends this framework to include $n$ variables represented as $x^i$.
\begin{example}
Consider a triple dynamical variable $\vec{u}=(x,y,z)$ that spans a three-dimensional phase space. Let $H$ and $G$ be Hamiltonians that determine the motion of points in the phase space. The fractal Nambu-Hamiltonians are given by the following equations
\begin{align}
  \frac{d_{F}^{\alpha} x}{d_{F}^{\alpha}t}&=\frac{\partial^{\alpha}_{F}(H,G)}
{\partial^{\alpha}_{F}(y,z )},\nonumber\\
\frac{d_{F}^{\alpha} y}{d_{F}^{\alpha}t}&=\frac{\partial^{\alpha}_{F}(H,G)}
{\partial^{\alpha}_{F}(z,x )},\nonumber\\
\frac{d_{F}^{\alpha} z}{d_{F}^{\alpha}t}&=\frac{\partial^{\alpha}_{F}(H,G)}
{\partial^{\alpha}_{F}(x,y )}.
\end{align}
Alternatively, in vector notation, we have:
\begin{equation}\label{ooppppij}
  \frac{d_{F}^{\alpha} \vec{u}}{d_{F}^{\alpha}t}=\vec{\nabla}^{\alpha} H\times \vec{\nabla}^{\alpha} G,
\end{equation}
where
\begin{equation}
  \vec{\nabla}^{\alpha}=D_{F,x}^{\alpha} \hat{e}_{i}+D_{F,y}^{\alpha} \hat{e}_{j}+D_{F,z}^{\alpha} \hat{e}_{k}.
\end{equation}
For any function $N(x,y,z)$, we have the following expression:
\begin{align}\label{uuhbvf}
  \frac{d_{F}^{\alpha} N}{d_{F}^{\alpha}t}&=\frac{\partial^{\alpha}_{F}(N,H,G)}
{\partial^{\alpha}_{F}(x,y,z )}\nonumber\\
&=\vec{\nabla}^{\alpha} N.(\vec{\nabla}^{\alpha} H\times \vec{\nabla}^{\alpha} G).
\end{align}
We can refer to the right-hand side of Eq. \eqref{uuhbvf} as the fractal Nambu-Poisson bracket, denoted by $[M,H,G]^{\alpha}$. The orbit of the system in the fractal phase space is determined by the intersection of two surfaces, namely $H=const.$ and $G=const.$. By using Eq. \eqref{ooppppij}, we can write:
\begin{equation}
 \vec{\nabla}^{\alpha}.( \vec{\nabla}^{\alpha} H\times \vec{\nabla}^{\alpha} G)=0.
\end{equation}
This may called the fractal Liouville theorem in the fractal Nambu phase space.
\end{example}
\begin{definition}
By selecting an arbitrary $H_r$ as a generalized Hamiltonian, we can define $J_{ij}$ as follows:
\begin{equation}
  J_{ij}=\epsilon_{ijk...l} D_{F,x_{k}}^{\alpha}H_{1}...D_{F,x_{z}}^{\alpha}H_{r-1}
D_{F,x_{c}}^{\alpha}H_{r+1}...D_{F,x_{l}}^{\alpha}H_{n-1},~~~H=H_{r}.
\end{equation}
In this definition, $J_{ij}$ satisfies the following conditions:
\begin{align}
  J_{ij} &= -J_{ji}, \\
  J_{il}D_{F,x_{l}}^{\alpha} J_{jk}+J_{jl}D_{F,x_{l}}^{\alpha} J_{ki}+J_{kl}D_{F,x_{l}}^{\alpha} J_{ij}&= 0,~~~(i,j,k,l=1,2,...,n).
\end{align}
\end{definition}
\begin{definition}
A mapping from $(x, y, z)$ to $(x', y', z')$ is considered a fractal canonical transform if it satisfies the condition:
\begin{equation}
  [x',y',z']=\frac{\partial^{\alpha}_{F}(x',y',z')}
{\partial^{\alpha}_{F}(x,y,z )}=1.
\end{equation}
Alternatively, it can be expressed as:
\begin{equation}
  \frac{\partial^{\alpha}_{F}(H',G')}
{\partial^{\alpha}_{F}(H,G)}=1,
\end{equation}
where $H' = h_{1}(H, G)$ and $G' = h_{2}(H, G)$.
In this definition, a fractal canonical transform refers to a mapping that preserves the fractal Nambu structure, ensuring that the fractal Nambu-Poisson brackets remain invariant under the transformation.
\end{definition}
\section{Fractal differential form and fractal Nambu mechanics  \label{S-3}}
In this section, we introduce the concept of fractal differential forms to derive Nambu's equations.
\begin{definition}
The fractal Poincar\'{e}-Cartan 1-form is defined as follows:
\begin{align}
  \Omega^{1}=pd_{F}^{\alpha}q-H(p,q)d_{F}^{\alpha}t.
\end{align}
\end{definition}
\begin{definition}
The fractal Poincar\'{e}-Cartan 2-form is defined as follows:
\begin{equation}
  \Omega^{2}=pd_{F}^{\alpha}q\wedge d_{F}^{\alpha}r-H_{1}d_{F}^{\alpha}H_{2}\wedge d_{F}^{\alpha}t.
\end{equation}
In these definitions, the variables $p$, $q$, and $r$ represent a triplet of dynamical variables, and $H_{1}$ and $H_{2}$ denote two Hamiltonians. The symbol $\wedge$ represents the wedge product, a mathematical operation used in multivariable calculus.
These fractal differential forms provide a mathematical framework for describing the dynamics and relationships between the variables involved in the Nambu system.
\end{definition}
The Hamilton and Nambu equations can be derived by taking the exterior differentials of $\Omega^{1}$ and $\Omega^{2}$, respectively. The resulting equations are as follows:\\
For $\Omega^{1}$:
\begin{align}
  d_{F}^{\alpha}\Omega^{1}=\left(d_{F}^{\alpha}p+
\frac{\partial^{\alpha}_{F} H}{\partial^{\alpha}_{F} q}d_{F}^{\alpha}t \right)\wedge \bigg(d_{F}^{\alpha}q -
\frac{\partial^{\alpha}_{F} H}{\partial^{\alpha}_{F} p}d_{F}^{\alpha}t\bigg)=0.
\end{align}
This equation leads to Hamilton's equations:
\begin{equation}
 \frac{d_{F}^{\alpha}q}{d_{F}^{\alpha}t}=\frac{\partial^{\alpha}_{F} H}{\partial^{\alpha}_{F} p},~~~\frac{d_{F}^{\alpha}p}{d_{F}^{\alpha}t}=-\frac{\partial^{\alpha}_{F} H}{\partial^{\alpha}_{F} q}.
\end{equation}
For $\Omega^{2}$:
\begin{align}
 & d_{F}^{\alpha}\Omega^{2}=\left(d_{F}^{\alpha}q-
\frac{\partial^{\alpha}_{F} (H_{1},H_{2})}{\partial^{\alpha}_{F} (p,r)}d_{F}^{\alpha}t \right)\wedge \left(d_{F}^{\alpha}p-
\frac{\partial^{\alpha}_{F} (H_{1},H_{2})}{\partial^{\alpha}_{F} (r,q)}d_{F}^{\alpha}t \right)\nonumber\\&\wedge \left(d_{F}^{\alpha}r-\frac{\partial^{\alpha}_{F} (H_{1},H_{2})}{\partial^{\alpha}_{F} (q,p)}d_{F}^{\alpha}t \right)=0.
\end{align}
This equation leads to the Nambu's equations:
\begin{align}\label{dds258}
\frac{d_{F}^{\alpha}q}{d_{F}^{\alpha}t}&=
\frac{\partial^{\alpha}_{F} (H_{1},H_{2})}{\partial^{\alpha}_{F} (p,r)} \nonumber\\
\frac{d_{F}^{\alpha}p}{d_{F}^{\alpha}t}&=
\frac{\partial^{\alpha}_{F} (H_{1},H_{2})}{\partial^{\alpha}_{F} (r,q)} \nonumber\\
\frac{d_{F}^{\alpha}r}{d_{F}^{\alpha}t}&=
\frac{\partial^{\alpha}_{F} (H_{1},H_{2})}{\partial^{\alpha}_{F} (q,p)}.
\end{align}
These equations describe the dynamics of the system in terms of the fractal derivatives and partial derivatives with respect to the variables involved. The Hamilton equations correspond to the case when $H_{1}$ and $H_{2}$ are the Hamiltonians, while the Nambu equations generalize the dynamics to more complex systems involving multiple variables and Hamiltonians.
\section{Fractal Lagrangians in Nambu mechanics}
We suggest a postulate that states the existence of Lagrangians matching to each Hamiltonian in order to achieve genuine Lagrangians in Nambu mechanics. We assume that the system contains an equal number of Lagrangians and Hamiltonians. Furthermore, we presume that each Lagrangian is consistent with the least action principle. In the most straightforward case, we suppose the presence of two Lagrangians, designated as $L_{1}(q(t), D^{\alpha}_{F}q )$ and $L_{2}(q(t), D^{\alpha}_{F}q )$, which are equivalent to the system's two Hamiltonians. With the help of this postulation, we can provide a thorough framework for figuring out the Lagrangians in Nambu mechanics \cite{ogawa2000nambu}.
\begin{align}
  \delta \int L_{1}(q(t), D^{\alpha}_{F}q )d_{F}^{\alpha}t&=0 \nonumber\\
 \delta \int L_{2}(q(t), D^{\alpha}_{F}q )d_{F}^{\alpha}t&=0
\end{align}
Thus we have the two Euler-Lagrange equations for the two Lagrangians:
\begin{align}\label{hhhygtgcxz}
 \frac{\partial^{\alpha}_{F} L_{1}}{\partial^{\alpha}_{F} t}-\frac{d_{F}^{\alpha}}{d_{F}^{\alpha}t}\frac{\partial^{\alpha}_{F} L_{1}}{\partial^{\alpha}_{F} D^{\alpha}_{F}q}&=0\\
\frac{\partial^{\alpha}_{F} L_{2}}{\partial^{\alpha}_{F} t}-\frac{d_{F}^{\alpha}}{d_{F}^{\alpha}t}\frac{\partial^{\alpha}_{F} L_{2}}{\partial^{\alpha}_{F} D^{\alpha}_{F}q}&=0.
\end{align}
Next, we define the first canonical momentum, denoted as $p$, and the second canonical momentum, denoted as $r$. The expressions for these momenta in terms of the Lagrangians are as follows:
\begin{align}\label{uuhbvgt}
p &= \frac{\partial^{\alpha}_{F} L_{1}}{\partial^{\alpha}_{F} D^{\alpha}_{F}q}\\
r &= \frac{\partial^{\alpha}_{F} L_{2}}{\partial^{\alpha}_{F} D^{\alpha}_{F}q}.
\end{align}
These momenta represent important quantities in Nambu mechanics, providing a measure of the conjugate momenta associated with the respective Lagrangians.
Using Eq.\eqref{hhhygtgcxz} we have
\begin{align}\label{qqawwuuhbvgt}
D^{\alpha}_{F}p &= \frac{\partial^{\alpha}_{F} L_{1}}{\partial^{\alpha}_{F} q}\\
D^{\alpha}_{F}r &= \frac{\partial^{\alpha}_{F} L_{2}}{\partial^{\alpha}_{F} q}.
\end{align}
Now we define the two Hamiltonians $H_1$ and $H_2$ through fractal exterior differentials
as follows:
\begin{equation}\label{33frdd}
  d_{F}^{\alpha}H_{1}\wedge d_{F}^{\alpha}H_{2}=\frac{1}{D^{\alpha}_{F}q}d_{F}^{\alpha}
(pD^{\alpha}_{F}q-L_{1})\wedge d_{F}^{\alpha}
(rD^{\alpha}_{F}q-L_{2}),
\end{equation}
where $H_1$ and $H_2$ are expressed in terms of $q, p$, and $r$. Then the left-hand side of Eq.\eqref{33frdd} is equal
 \begin{align}\label{dttyyuuzz}
   d_{F}^{\alpha}H_{1}\wedge d_{F}^{\alpha}H_{2}&=\left(\frac{\partial^{\alpha}_{F} H_{1}}{\partial^{\alpha}_{F} q} d_{F}^{\alpha}q+\frac{\partial^{\alpha}_{F} H_{1}}{\partial^{\alpha}_{F} p} d_{F}^{\alpha}p+\frac{\partial^{\alpha}_{F} H_{1}}{\partial^{\alpha}_{F} r} d_{F}^{\alpha}r\right)\nonumber\\&\wedge \left(\frac{\partial^{\alpha}_{F} H_{2}}{\partial^{\alpha}_{F} q} d_{F}^{\alpha}q+\frac{\partial^{\alpha}_{F} H_{2}}{\partial^{\alpha}_{F} p} d_{F}^{\alpha}p+\frac{\partial^{\alpha}_{F} H_{2}}{\partial^{\alpha}_{F} r} d_{F}^{\alpha}r\right)\nonumber\\&=
\bigg(\frac{\partial^{\alpha}_{F} (H_{1},H_{2})}{\partial^{\alpha}_{F} (p,r)} d_{F}^{\alpha}p\wedge d_{F}^{\alpha}r+\frac{\partial^{\alpha}_{F} (H_{1},H_{2})}{\partial^{\alpha}_{F} (r,q)} d_{F}^{\alpha}r\wedge d_{F}^{\alpha}q \nonumber\\&+\frac{\partial^{\alpha}_{F} (H_{1},H_{2})}{\partial^{\alpha}_{F} (q,p)} d_{F}^{\alpha}q\wedge d_{F}^{\alpha}p\bigg),
 \end{align}
and the right-hand side of Eq.\eqref{33frdd} is equal to, by virtue of Eqs. \eqref{uuhbvgt} and
Eqs. \eqref{qqawwuuhbvgt} ,
\begin{align}\label{swwwakkss}
&  \frac{1}{D^{\alpha}_{F}q}d_{F}^{\alpha}
(pD^{\alpha}_{F}q-L_{1})\wedge d_{F}^{\alpha}
(rD^{\alpha}_{F}q-L_{2})\nonumber\\&=\frac{1}{D^{\alpha}_{F}q}
\left[D^{\alpha}_{F}q d_{F}^{\alpha}p+
pd_{F}^{\alpha}D^{\alpha}_{F}q-\left(\frac{\partial^{\alpha}_{F} L_{1}}{\partial^{\alpha}_{F} q}d_{F}^{\alpha}q+ \frac{\partial^{\alpha}_{F} L_{1}}{\partial^{\alpha}_{F} D^{\alpha}_{F}q}d_{F}^{\alpha}D^{\alpha}_{F}q\right)\right]\nonumber\\&\wedge
\left[D^{\alpha}_{F}q d_{F}^{\alpha}r+
rd_{F}^{\alpha}D^{\alpha}_{F}q-\left(\frac{\partial^{\alpha}_{F} L_{2}}{\partial^{\alpha}_{F} q}d_{F}^{\alpha}q+ \frac{\partial^{\alpha}_{F} L_{2}}{\partial^{\alpha}_{F} D^{\alpha}_{F}q}d_{F}^{\alpha}D^{\alpha}_{F}q\right)\right]\nonumber\\&=
D^{\alpha}_{F}qd_{F}^{\alpha}p\wedge d_{F}^{\alpha}r+D^{\alpha}_{F}pd_{F}^{\alpha}r\wedge d_{F}^{\alpha}q+
D^{\alpha}_{F}rd_{F}^{\alpha}q\wedge d_{F}^{\alpha}p.
\end{align}
Comparing Eqs.\eqref{dttyyuuzz} and \eqref{swwwakkss}, we obtain Eq.\eqref{dds258} that establish the validity of our postulate regarding the existence of two Lagrangians and the corresponding definitions of the Hamiltonians.
\section{Applications of Fractal Nambu Mechanics \label{S-4}}
In this section, we give some applications in classical mechanics.
\begin{example}
Consider the Euler asymmetric top as an example. The Euler equation for this system can be derived as follows. Let's define the following Hamiltonians:
\begin{equation}
  H_{1}=\frac{L_{1}^{2}}{2I_{1}}+
\frac{L_{2}^{2}}{2I_{2}}+\frac{L_{3}^{2}}{2I_{3}},~~~H_{2}=
\frac{1}{2}(L_{1}^{2}+L_{2}^{2}+L_{3}^{2}),
\end{equation}
where $I_{i},~i=1,2,3$ represent the inertia, and $L_{i},~i=1,2,3$ represent the components of the angular momenta. Using Eq.\eqref{uuhbvf}, we can obtain fractal Euler equations as:
\begin{equation}
 \frac{d_{F}^{\alpha} L_{i}}{d_{F}^{\alpha}t}=\frac{\partial^{\alpha}_{F}(H_{1},H_{2})}
{\partial^{\alpha}_{F}(x,y,z )}=\epsilon_{ijk}(D_{F,x_{j}}^{\alpha}H_{1}D_{F,x_{k}}^{\alpha}H_{2}).
\end{equation}
Alternatively, we can express it as:
\begin{equation}
 \frac{d_{F}^{\alpha} L_{1}}{d_{F}^{\alpha}t}=
\frac{I_{3}-I_{2}}{I_{2}I_{3}}L_{2}L_{3},~~~\frac{d_{F}^{\alpha} L_{2}}{d_{F}^{\alpha}t}=
\frac{I_{1}-I_{3}}{I_{1}I_{3}}L_{2}L_{3},~~~\frac{d_{F}^{\alpha} L_{3}}{d_{F}^{\alpha}t}=
\frac{I_{2}-I_{1}}{I_{1}I_{2}}L_{1}L_{2}.
\end{equation}
These equations describe the fractal Euler equations for the Euler asymmetric top system, where the fractal time derivatives of the angular momenta $L_{i}$ are expressed in terms of the inertia $I_{i}$ and the angular momenta themselves.
\end{example}
\begin{example}
The fractal Nahm system, which serves as a model for static monopoles, can be obtained using fractal Nambu's equation. Consider the Hamiltonians $H_{1}=a_{1}x_{1}^2-a_{2}x_{2}^2$ and $H_{2}=a_{1}x_{1}^2-a_{3}x_{3}^2$ where $a_{1},~a_{2},~a_{3}$ are constant. Applying fractal Nambu's equation, we have:
\begin{equation}
  \frac{d_{F}^{\alpha} x_{i}}{d_{F}^{\alpha}t}=\{H_{1},H_{2},x_{i}\}.
\end{equation}
Alternatively, we can express it as:
\begin{equation}\label{rrrfff}
  \frac{d_{F}^{\alpha} x_{1}}{d_{F}^{\alpha}t}=a_{2}a_{3}x_{2}x_{3},~~~\frac{d_{F}^{\alpha} x_{2}}{d_{F}^{\alpha}t}=a_{1}a_{3}x_{1}x_{3},~~~\frac{d_{F}^{\alpha} x_{3}}{d_{F}^{\alpha}t}=a_{1}a_{2}x_{1}x_{2}.
\end{equation}

\begin{figure}
        \includegraphics[width=\textwidth]{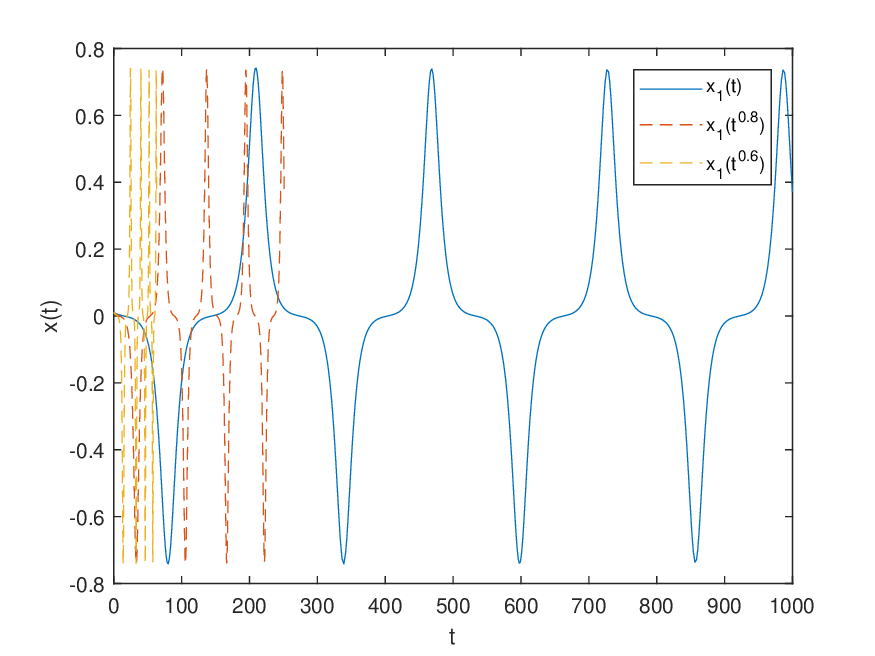}
        \caption{Plot of $x_{1}(t)$ exhibiting various fractal dimensions.}
        \label{uu1}
        \end{figure}
 \begin{figure}
        \includegraphics[width=\textwidth]{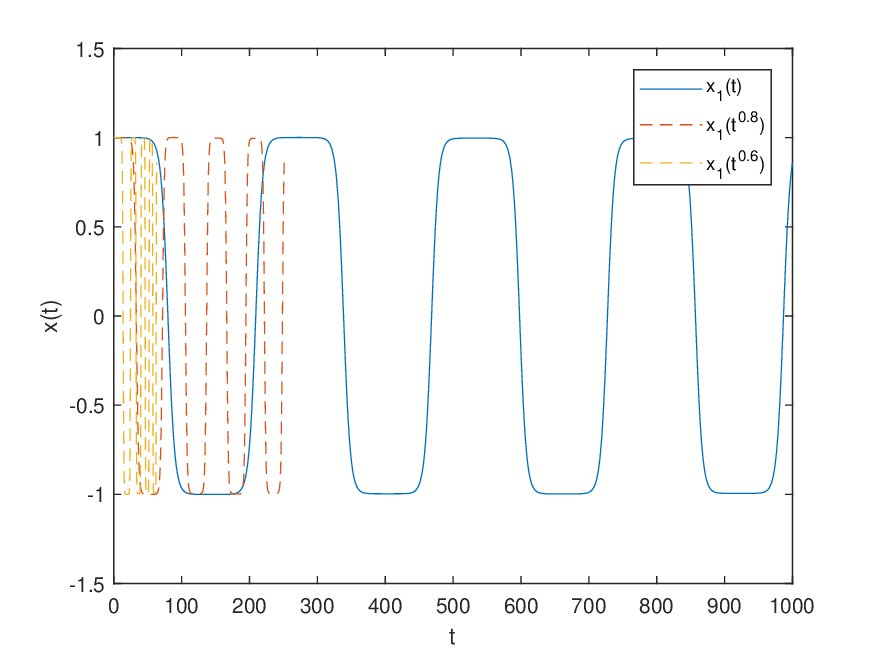}
        \caption{Plot of $x_{2}(t)$ exhibiting various fractal dimensions}
        \label{uu2}
 \end{figure}
\begin{figure}
        \includegraphics[width=\textwidth]{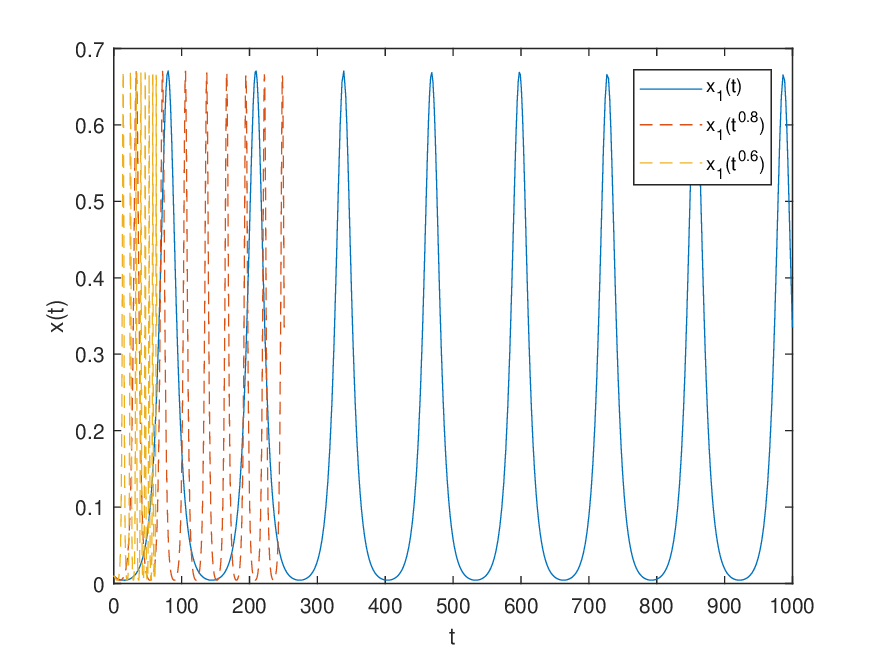}
        \caption{Plot of $x_{3}(t)$ for various fractal dimensions}
        \label{uu3}
\end{figure}

In Figures \ref{uu1}, \ref{uu2}, and \ref{uu3}, we showcase the solutions to Eq.\eqref{rrrfff}. Here, we apply the condition $S_{F}^{\alpha}(t)\leq t^{\alpha}$ with specific values assigned to the parameters, namely $a_{1}=0.40452$, $a_{2}= -0.222486$, and $a_{3}=0.494413$.

These equations describe the fractal Nahm system, providing the time derivatives of the variables $x_{i}$ as products of the corresponding variables themselves. This model captures the dynamics of static monopoles.
\end{example}
\begin{example}
Consider a 2-dimensional fractal harmonic oscillator described by the following Hamiltonians:
\begin{align}
H_{1}&=p_{1}^{2}+x_{1}^2=C_{1}\nonumber\\
H_{2}&=p_{2}^{2}+x_{2}^2=C_{2}\nonumber\\
H_{3}&=x_{1}p_{2}-x_{2}p_{1}=C_{3}\nonumber\\
H_{4}&=p_{1}p_{2}+x_{1}x_{2}=C_{4}.
\end{align}
Here, $C_{i}, i=1,2,3,4$, are constants. The corresponding fractal Nambu's equation is given by:
\begin{equation}
  \{H_{1},H_{2},H_{3},f\}=-\frac{1}{2C_{4}}
\frac{\partial^{\alpha}_{F}(H_{1},H_{2},H_{2},H_{3},f)}
{\partial^{\alpha}_{F}(p_{1},p_{2},x_{1},x_{2} )}.
\end{equation}
In this equation, the curly brackets represent the fractal Nambu-Poisson bracket, and $f$ is an arbitrary function. The equation captures the dynamics of the system and highlights the interplay between the Hamiltonians and the additional function $f$. The fractal Nambu's equation provides a powerful tool for studying the behavior of the 2-dimensional fractal harmonic oscillator system.
\end{example}
\section{Conclusion  \label{S-5}}
In conclusion, Fractal Nambu mechanics is a mathematical framework that extends Nambu mechanics to incorporate fractal geometries and dynamics. Fractals, which are shapes with self-similar properties and fractional dimensions, are found in various natural and artificial systems. By utilizing fractal derivatives, fractal differential forms, and the fractal Nambu-Poisson bracket, Fractal Nambu mechanics provides a powerful tool for describing the dynamics of complex systems with fractal characteristics.
The key features of Fractal Nambu mechanics include the formulation of fractal Nambu equations, which generalize the Hamilton and Nambu equations to systems with multiple variables and Hamiltonians. The fractal Nambu-Poisson bracket satisfies properties such as skew-symmetry, the Jacobi identity, and the Leibniz rule, ensuring consistency and preserving the Nambu structure. Fractal canonical transforms, which preserve the fractal Nambu structure, allow for the study of different coordinate systems.
The introduction of fractal differential forms, such as the fractal Poincar\'{e}-Cartan 1-form and 2-form, enables the derivation of Nambu's equations and provides a mathematical framework for describing the dynamics and relationships between variables in the Nambu system.
Fractal Nambu mechanics finds applications in various areas of classical mechanics, allowing for the study of complex systems with fractal properties. Examples include the Euler asymmetric top, where the fractal Euler equations can be derived, and the study of dynamical systems in fractal phase spaces.
Overall, Fractal Nambu mechanics provides a comprehensive and powerful framework for understanding and modeling complex systems with fractal characteristics, offering new insights and approaches to studying the dynamics of fractal structures.

\bibliographystyle{plain}

\bibliography{fractalnumbu}

\begin{thebibliography}{10}

\bibitem{agathiyan2023integral}
A~Agathiyan, A~Gowrisankar, and Nur Aisyah~Abdul Fataf.
\newblock On the integral transform of fractal interpolation functions.
\newblock {\em Mathematics and Computers in Simulation}, 2023.

\bibitem{agathiyan2023remarks}
A~Agathiyan, A~Gowrisankar, Nur Aisyah~Abdul Fataf, and Jinde Cao.
\newblock Remarks on the integral transform of non-linear fractal interpolation functions.
\newblock {\em Chaos, Solitons \& Fractals}, 173:113749, 2023.

\bibitem{ma-5}
Alexander~S. Balankin.
\newblock A continuum framework for mechanics of fractal materials i: from fractional space to continuum with fractal metric.
\newblock {\em Eur. Phys. J. B}, 88(4), April 2015.

\bibitem{balankin2023vector}
Alexander~S Balankin and Baltasar Mena.
\newblock Vector differential operators in a fractional dimensional space, on fractals, and in fractal continua.
\newblock {\em Chaos Solit. Fractals}, 168:113203, 2023.

\bibitem{baleanu2009fractional}
Dumitru Baleanu, Alireza~K Golmankhaneh, and Ali~K Golmankhaneh.
\newblock Fractional nambu mechanics.
\newblock {\em International Journal of Theoretical Physics}, 48:1044--1052, 2009.

\bibitem{banchuin20224noise}
Rawid Banchuin.
\newblock Nonlocal fractal calculus based analyses of electrical circuits on fractal set.
\newblock {\em COMPEL - Int. J. Comput. Math. Electr. Electron. Eng.}, 41(1):528--549, 2022.

\bibitem{ma-7}
Martin~T. Barlow and Edwin~A. Perkins.
\newblock Brownian motion on the sierpinski gasket.
\newblock {\em Probab. Theory Rel.}, 79(4):543--623, nov 1988.

\bibitem{Barnsley-book}
Michael~F Barnsley.
\newblock {\em Fractals Everywhere}.
\newblock Academic Press, New York, 2014.

\bibitem{bayen1975remarks}
F~Bayen and M~Flato.
\newblock Remarks concerning nambu's generalized mechanics.
\newblock {\em Physical Review D}, 11(10):3049, 1975.

\bibitem{bevilacqua2023inverse}
Luiz Bevilacqua and Marcelo~M Barros.
\newblock The inverse problem for fractal curves solved with the dynamical approach method.
\newblock {\em Chaos Solit. Fractals}, 168:113113, 2023.

\bibitem{bishop2017fractals}
Christopher~J Bishop and Yuval Peres.
\newblock {\em Fractals in probability and analysis}, volume 162.
\newblock Cambridge University Press, Cambridge, 2017.

\bibitem{bongiorno2011henstock}
D~Bongiorno, G~Corrao, et~al.
\newblock The {H}enstock-{K}urzweil-{S}tieltjes type integral for real functions on a fractal subset of the real line.
\newblock In {\em Bollettino di Matematica pura e applicata vol. IV}, volume~4, pages 5--16. Aracne, Boston, 2011.

\bibitem{bongiorno2018derivatives}
Donatella Bongiorno.
\newblock Derivatives not first return integrable on a fractal set.
\newblock {\em Ricerche di Matematica}, 67(2):597--604, 2018.

\bibitem{Bongiorno23}
Donatella Bongiorno.
\newblock Derivation and integration on a fractal subset of the real line.
\newblock In Dr. Sid-Ali Ouadfeul, editor, {\em Fractal Analysis - Applications and Updates}, chapter~7. IntechOpen, Rijeka, 2023.

\bibitem{bongiorno2015fundamental}
Donatella Bongiorno and Giuseppa Corrao.
\newblock On the fundamental theorem of calculus for fractal sets.
\newblock {\em Fractals}, 23(02):1550008, 2015.

\bibitem{chandre2023classical}
Cristel Chandre and Atsushi Horikoshi.
\newblock Classical nambu brackets in higher dimensions.
\newblock {\em J. Math. Phys.}, 64(5), 2023.

\bibitem{cohen1975nambu}
Isaac Cohen and Andr{\'e}s~J K{\'a}lnay.
\newblock On nambu's generalized hamiltonian mechanics.
\newblock {\em International Journal of Theoretical Physics}, 12:61--67, 1975.

\bibitem{ma-9}
M.~Czachor.
\newblock Waves along fractal coastlines: From fractal arithmetic to wave equations.
\newblock {\em Acta Phys. Pol. B}, 50(4):813, 2019.

\bibitem{czachor2023contra}
M~Czachor.
\newblock Contra bellum: Bell's theorem as a confusion of languages.
\newblock {\em Acta Phys. Pol}, 143(6):S158--S158, 2023.

\bibitem{damian2023mechanical}
Lucero Dami{\'a}n~Adame, Claudia del~Carmen Guti{\'e}rrez-Torres, Bernardo Figueroa-Espinoza, Juan~Gabriel Barbosa-Salda{\~n}a, and Jos{\'e}~Alfredo Jim{\'e}nez-Bernal.
\newblock A mechanical picture of fractal darcy’s law.
\newblock {\em Fractal and Fractional}, 7(9):639, 2023.

\bibitem{deppman2022tsallis}
Airton Deppman and Eugenio Meg{\'\i}as.
\newblock Tsallis statistics and {QCD} thermodynamics.
\newblock In {\em EPJ Web of Conferences}, volume 270, page 00033, 2022.

\bibitem{el1993certain}
MS~El~Naschie.
\newblock On certain infinite dimensional cantor sets and the schr{\"o}dinger wave.
\newblock {\em Chaos Solit. Fractals}, 3(1):89--98, 1993.

\bibitem{freiberg2002harmonic}
U~Freiberg and Martina Z{\"a}hle.
\newblock Harmonic calculus on fractals-a measure geometric approach {I}.
\newblock {\em Potential Anal.}, 16(3):265--277, 2002.

\bibitem{gautheron1996some}
Philippe Gautheron.
\newblock Some remarks concerning nambu mechanics.
\newblock {\em Letters in Mathematical Physics}, 37:103--116, 1996.

\bibitem{giona1995fractal}
Massimiliano Giona.
\newblock Fractal calculus on [0, 1].
\newblock {\em Chaos Solit. Fractals}, 5(6):987--1000, 1995.

\bibitem{khalili2011fractional}
Ali~Khalili Golmankhaneh, AlirezaKhalili Golmankhaneh, Dumitru Baleanu, and Mihaela~Cristina Baleanu.
\newblock Fractional odd-dimensional mechanics.
\newblock {\em Adv. Differ. Equ.}, 2011:1--12, 2011.

\bibitem{Alireza-book}
Alireza~Khalili Golmankhaneh.
\newblock {\em Fractal Calculus and its Applications}.
\newblock World Scientific, Singapore, 2022.

\bibitem{khalili2023non}
Alireza~Khalili Golmankhaneh, Kerri Welch, Cristina Serpa, and Palle E~T J{\o}rgensen.
\newblock Non-standard analysis for fractal calculus.
\newblock {\em The Journal of Analysis}, 31:1895–1916, 2023.

\bibitem{golmankhaneh2023classical}
Alireza~Khalili Golmankhaneh, Kerri Welch, Cemil Tun{\c{c}}, and Yusif~S Gasimov.
\newblock Classical mechanics on fractal curves.
\newblock {\em Eur. Phys. J. Spec. Top.}, pages 1--9, 2023.

\bibitem{gowrisankar2021fractal}
Arulprakash Gowrisankar, Alireza~Khalili Golmankhaneh, and Cristina Serpa.
\newblock Fractal calculus on fractal interpolation functions.
\newblock {\em Fractal Fract.}, 5(4):157, 2021.

\bibitem{grabowski1999remarks}
Janusz Grabowski and Giuseppe Marmo.
\newblock Remarks on nambu-poisson and nambu-jacobi brackets.
\newblock {\em J. Phys. A Math.}, 32(23):4239, 1999.

\bibitem{guha2002applications}
Partha Guha.
\newblock Applications of nambu mechanics to systems of hydrodynamical type.
\newblock {\em J. Math. Phys.}, 43(8):4035--4040, 2002.

\bibitem{guha2004applications}
Partha Guha.
\newblock Applications of nambu mechanics to systems of hydrodynamical type ii.
\newblock {\em J. Nonlinear Math. Phys.}, 11(2):223--232, 2004.

\bibitem{hambly2023dimension}
Ben Hambly and Peter Koepernik.
\newblock Dimension results and local times for superdiffusions on fractals.
\newblock {\em Stoch. Process. their Appl.}, 158:377--417, 2023.

\bibitem{heider2022self}
Yousef Heider, Franz Bamer, Firaz Ebrahem, and Bernd Markert.
\newblock Self-organized criticality in fracture models at different scales.
\newblock {\em Examples and Counterexamples}, 2:100054, 2022.

\bibitem{jiang1998some}
Huikun Jiang and Weiyi Su.
\newblock Some fundamental results of calculus on fractal sets.
\newblock {\em Commun. Nonlinear Sci. Numer. Simul.}, 3(1):22--26, 1998.

\bibitem{ma-8}
Jun Kigami.
\newblock {\em Analysis on Fractals}.
\newblock Cambridge University Press, Cambridge, jun 2001.

\bibitem{kotov2023application}
MA~Kotov, VM~Volchuk, DM~Zeziukov, and TM~Pavlenko.
\newblock Application of 3d fractal modeling for predicting concrete strength characteristics.
\newblock In {\em AIP Conference Proceedings}, volume 2526. AIP Publishing, 2023.

\bibitem{Michel-j}
Michel~L. Lapidus, Goran Radunovi{\'{c}}, and Darko {\v{Z}}ubrini{\'{c}}.
\newblock {\em Fractal Zeta Functions and Fractal Drums}.
\newblock Springer International Publishing, New York, 2017.

\bibitem{lassig1997constrained}
CC~Lassig and Girish~C Joshi.
\newblock Constrained systems described by nambu mechanics.
\newblock {\em Lett. Math. Phys .}, 41:59--63, 1997.

\bibitem{linton2021fractals}
O.~Linton.
\newblock {\em Fractals: On the Edge of Chaos}.
\newblock Wooden Books. Bloomsbury, New York, 2021.

\bibitem{llibre2023generalized}
Jaume Llibre, Rafael Ram{\'\i}rez, and Valent{\'\i}n Ram{\'\i}rez.
\newblock Generalized {C}artesian--{N}ambu vector fields.
\newblock In {\em Dynamics through First-Order Differential Equations in the Configuration Space}, pages 177--283. Springer, New York, 2023.

\bibitem{makhaldiani2011nambu}
Nugzar Makhaldiani.
\newblock Nambu-poisson dynamics with some applications.
\newblock {\em arXiv preprint arXiv:1201.5105}, 2011.

\bibitem{Mandelbrot-1}
Benoit~B Mandelbrot.
\newblock {\em The Fractal Geometry of Nature}.
\newblock WH freeman, New York, 1982.

\bibitem{megias2022nambu}
E~Megias, MJ~Teixeira, VS~Tim{\'o}teo, and A~Deppman.
\newblock Nambu--jona-lasinio model with a fractal inspired coupling.
\newblock {\em arXiv preprint arXiv:2203.11080}, 2022.

\bibitem{nambu1952lagrangian}
Yoichiro Nambu.
\newblock On lagrangian and hamiltonian formalism.
\newblock {\em Progress of Theoretical Physics}, 7(2):131--170, 1952.

\bibitem{nambu1973generalized}
Yoichiro Nambu.
\newblock Generalized hamiltonian dynamics.
\newblock {\em Physical Review D}, 7(8):2405, 1973.

\bibitem{Nottale-4}
L.~Nottale and J.~Schneider.
\newblock Fractals and nonstandard analysis.
\newblock {\em J. Math. Phys.}, 25(5):1296--1300, may 1984.

\bibitem{nottale1993fractal}
Laurent Nottale.
\newblock {\em Fractal Space-Time and Microphysics: Towards a Theory of Scale Relativity}.
\newblock World Scientific, Singapore, 1993.

\bibitem{ogawa2000nambu}
Tsuyoshi Ogawa and Toshiaki Sagae.
\newblock Nambu mechanics in the lagrangian formalism.
\newblock {\em International Journal of Theoretical Physics}, 39:2875--2890, 2000.

\bibitem{ord1983fractal}
GN~Ord.
\newblock Fractal space-time: a geometric analogue of relativistic quantum mechanics.
\newblock {\em J. Phys. A Math.}, 16(9):1869, 1983.

\bibitem{pandit1998generalized}
Sagar~A Pandit and Anil~D Gangal.
\newblock On generalized nambu mechanics.
\newblock {\em J. Phys. A Math.}, 31(12):2899, 1998.

\bibitem{AD-2}
A.~Parvate and A.D. Gangal.
\newblock Calculus on fractal subsets of real line-{II}: Conjugacy with ordinary calculus.
\newblock {\em Fractals}, 19(03):271--290, sep 2011.

\bibitem{parvate2011calculus}
A.~Parvate, S.~Satin, and A.D. Gangal.
\newblock Calculus on fractal curves in $\mathbb{R}^{n}$.
\newblock {\em Fractals}, 19(01):15--27, mar 2011.

\bibitem{parvate2009calculus}
Abhay Parvate and Anil~D Gangal.
\newblock Calculus on fractal subsets of real line-{I}: Formulation.
\newblock {\em Fractals}, 17(01):53--81, 2009.

\bibitem{priyanka2023alpha}
TMC Priyanka, C~Serpa, and A~Gowrisankar.
\newblock $\alpha$-fractal function with variable parameters: An explicit representation.
\newblock {\em Fractals}, 2023.

\bibitem{samayoa2022hausdorff}
Didier Samayoa, Ernesto Pineda~Le{\'o}n, Lucero Dami{\'a}n~Adame, Eduardo Reyes~de Luna, and Andriy Kryvko.
\newblock The hausdorff dimension and capillary imbibition.
\newblock {\em Fractal and Fractional}, 6(6):332, 2022.

\bibitem{trifce2020fractional}
Trifce. Sandev.
\newblock {\em Fractional Equations and Models: Theory and Applications}.
\newblock Springer, New York, 2020.

\bibitem{shishanin2018nambu}
AO~Shishanin.
\newblock Nambu mechanics and its applications.
\newblock In {\em IOP Conference Series: Materials Science and Engineering}, volume 468, page 012029. IOP Publishing, 2018.

\bibitem{Shlesinger-6}
M~F Shlesinger.
\newblock Fractal time in condensed matter.
\newblock {\em Annu. Rev. Phys. Chem.}, 39(1):269--290, oct 1988.

\bibitem{stillinger1977axiomatic}
Frank~H Stillinger.
\newblock Axiomatic basis for spaces with noninteger dimension.
\newblock {\em J. Math. Phys.}, 18(6):1224--1234, 1977.

\bibitem{takhtajan1994foundation}
Leon Takhtajan.
\newblock On foundation of the generalized nambu mechanics.
\newblock {\em Commun. Math. Phys.}, 160:295--315, 1994.

\bibitem{ma-6}
Vasily~E. Tarasov.
\newblock {\em Fractional Dynamics}.
\newblock Springer Berlin Heidelberg, New York, 2010.

\bibitem{tarraf2023fractal}
Walid Tarraf, Diogo Queiros-Cond{\'e}, Patrick Ribeiro, and Rafik Absi.
\newblock Fractal geometric model for statistical intermittency phenomenon.
\newblock {\em Entropy}, 25(5):749, 2023.

\bibitem{uchaikin2013fractional}
Vladimir~V Uchaikin.
\newblock {\em Fractional Derivatives for Physicists and Engineers}, volume~2.
\newblock Springer, New York, 2013.

\bibitem{valarmathi2023variable}
R~Valarmathi and A~Gowrisankar.
\newblock Variable order fractional calculus on $\alpha$-fractal functions.
\newblock {\em The Journal of Analysis}, 31:2799--2815, 2023.

\bibitem{Vrobel}
Susie Vrobel.
\newblock {\em Fractal Time}.
\newblock World Scientific, Singapore, jan 2011.

\bibitem{withers1988fundamental}
W.D. Withers.
\newblock Fundamental theorems of calculus for hausdorff measures on the real line.
\newblock {\em J. Math. Anal. Appl.}, 129(2):581--595, feb 1988.

\bibitem{xu2015fractional}
Yan-Li Xu and Shao-Kai Luo.
\newblock Fractional nambu dynamics.
\newblock {\em Acta Mechanica}, 226(11):3781--3793, 2015.

\end{thebibliography}

\end{document}